\documentclass[11pt]{article}
\usepackage{amsfonts}
\begin{document}

\title{\Large
Integrable geodesic flows on the suspensions of toric automorphisms
\thanks{Submitted to Proceedings of the Steklov Institute of Mathematics
(in Russian).}}
\author{\large Alexey V. BOLSINOV
\thanks{
Department of Mathematics and Mechanics,
    Moscow State University, 119899 Moscow, Russia, e-mail:
    bols@difgeo.math.msu.su}
\ and Iskander A. TAIMANOV
\thanks{Institute of Mathematics,
    630090 Novosibirsk, Russia,
    e-mail: taimanov@math.nsc.ru}
}
\date{}
\maketitle

\newtheorem{lemma}{Lemma}
\newtheorem{theorem}{Theorem}
\newtheorem{definition}{Definition}
\newtheorem{corollary}{Corollary}
\newtheorem{remark}{Remark}
\newtheorem{problem}{Problem}
\newtheorem{conjecture}{Conjecture}
\newtheorem{proposition}{Proposition}

\def\R{{\mathbb R}}
\def\C{{\mathbb C}}
\def\Z{{\mathbb Z}}
\def\H{{\mathbb H}}
\def\N{{\mathbb N}}
\def\Tr{{\rm Tr}\,}
\def\ad{{\rm ad}\,}
\def\Ad{{\rm Ad}\,}
\newcommand{\Ker}{{\rm Ker}\,}
\newcommand{\Sp}{{\rm Sp}\,}
\renewcommand{\Re}{{\rm Re}\,}
\renewcommand{\Im}{{\rm Im}\,}

\section{Introduction and main results}

\medskip

In this paper we resume our study of integrable geodesic flows on the
suspensions of toric automorphisms which we started in \cite{BT}.

A closed manifold $M_A = M^{n+1}$ is called the suspension of a toric
automorphism $A: T^n \to T^n$ if there is a fibration
\begin{equation}
\pi: M^{n+1} \to S^1
\label{fibration}
\end{equation}
of this manifold over the circle $S^1$ with $T^n$-fibres
such that the monodromy of this fibration is given by $A \in SL(n,\Z)$.

The manifold $M_A$ is constructed as the quotient of the free
$\Z$-action
$$
(X,z) \to (AX,z+1)
$$
on the cylinder $T^n \times \R$ where $X \in T^n = \R^n/\Z^n, z \in \R$.

\begin{theorem}
If all eigenvalues of an automorphism $A \in SL(n,\Z)$ are real or $n=2$,
then $M_A$, the suspension of $A$,
admits a real-analytic Riemannian metric such that

1) the geodesic flow of this metric is (Liouville) integrable in terms of
$C^{\infty}$ first integrals;

2) the measure entropy of the geodesic flow with respect to any smooth
invariant measure vanishes;

3) the topological entropy of this flow meets the following inequality
\begin{equation}
h_{\rm top} \geq \log \left(
\max_{\lambda \in \Sp A} |\lambda| \right),
\label{entropy}
\end{equation}
where $\Sp A$ is the spectrum of $A$, i.e. the set of its eigenvalues.
\end{theorem}

For $A$ the identity, $M_A$ is a torus and in this case the statement of
the theorem is evident.

The first nontrivial case was found by Butler \cite{Butler}
who constructed an integrable geodesic flow on the manifold $M_A$ with
$$
A =
\left(
\begin{array}{cc}
1 & 1 \\
0 & 1
\end{array}
\right).
$$
He constructed the metric as a homogeneous metric on a nilmanifold
and worked in terms of global coordinates on the corresponding nilpotent
Lie group. In particular, Butler
showed that some topological obstructions to integrability
of geodesic flows in terms of real-analytic (or in some sense geometrically
simple) first integrals found in \cite{T1,T2} do not obstruct integrability
in terms of $C^{\infty}$ functions.

The suspension construction was found in \cite{BT}. In this paper
generalizing Butler's analytic trick for constructing $C^{\infty}$
first integrals we constructed an integrable geodesic flow on the
manifold $M_A$ with
\begin{equation}
A =
\left(
\begin{array}{cc}
2 & 1 \\
1 & 1
\end{array}
\right)
\label{matrix}
\end{equation}
and had explained that this suspension construction is quite general.
In \cite{BT} we discussed only one concrete example which appears to be the
first example of Liouville integrable geodesic flow with positive
topological entropy and also the first example of the geodesic flow for which
the Liouville entropy vanishes but the topological entropy is positive.

We shall study the Lyapunov exponents of
the flow from \cite{BT} and prove the following
statement.

\begin{theorem}
Given the Riemannian manifold $M_A$ with $A$ of the form (\ref{matrix}) and
the metric constructed in \cite{BT} (see Section 4), the unit cotangent bundle
$S M_A$ contains two four-dimensional invariant submanifolds $N^u$ and
$N^v$ such that

1) $N^u$ and $N^v$ are diffeomorphic to $M_A \times S^1$;

2)  the intersection $N^u \cap N^v = V$ consists of two three-dimensional
components $V^+$ and $V^-$, each of these components is diffeomorphic to
$M_A$ and consists in trajectories orthogonal to the fibers of the fibration
(\ref{fibration});

3) the Lyapunov exponents vanish at points from
$S M_A \setminus \{N^u \cup N^v\}$ and for any point from
$N^u \cup N^v$ there are
nonzero Lyapunov exponents;

4) all invariant (Borel) measures on $N^u$ and $N^v$ are supported on
$V^+ \cup V^-$ and there are smooth invariant measures on $V^+$ and $V^-$;

5) $N^u$ is a stable manifold for $V^+$ and an unstable manifold
for $V^-$, i.e., any trajectory in $N^u \setminus V$ is asymptotic to
a trajectory from $V^+$ as $t \to \infty$ and is asymptotic to a trajectory
from $V^-$ as $t \to -\infty$;

6) $N^v$ is a stable manifold for $V^-$ and an unstable manifold for
$V^+$;

7) the complement to $N^u \cup N^v$ is fibered by invariant tori.
\end{theorem}

Now we derive from this theorem that

\begin{corollary}
Given the Riemannian manifold $M_A$ with $A$ of the form (\ref{matrix}) and
the metric constructed in \cite{BT} (see Section 4), the topological entropy
of this flow equals
$$
h_{\rm top} = \log \frac{3+\sqrt{5}}{2}
$$
and there are measures of maximal entropy supported on
$V^+$ or $V^-$.
\end{corollary}

We would like to mention the following property of this integrable flow:

\begin{corollary}
The restrictions of the geodesic flow on $M_A$
onto $V^+$ or $V^-$ are Anosov flows.
\end{corollary}

One can see that easily: take a fiber of the fibration (\ref{fibration})
and at each point $q$ of the fiber take a covector $p=(p_u=p_v=0, p_z = 1)$.
Such points $(q,p)$ form a two-torus $T^2$ embedded into $S M_A$.
Then draw a geodesic in the direction of this covector. After the unit time
it will return back to this fiber and therefore we have a recurrence mapping
$$
T^2 \to T^2
$$
given by the hyperbolic matrix (\ref{matrix}).

\section{Entropy and integrability}

In this section we recall some well known
definitions and facts from the theory of dynamical systems.
For detailed explanation of different facts from this section we refer to
\cite{BF,KH,Sinai}.

A) {\sl Geodesic flows as Hamiltonian systems.}

Let $M^n$ be a Riemannian manifold with the metric $g_{ij}$.
Denote local coordinates on the cotangent bundle $T^\ast M^n$ as
$(x^1,\dots,x^n,p_1,\dots,p_n)$ where $(x_1,\dots,x^n)$ are
(local) coordinates on $M^n$ and the momenta $p_1,\dots,p_n$ are
defined from tangent vectors (velocities of curves on $M^n$) by
the Legendre transformation:
$$
p_i = g_{ij}\dot{x}^j.
$$
There is a symplectic form
$$
\omega = \sum_{i=1}^n d x^i \wedge d p_i
$$
on $T^\ast M^n$ which is correctly defined globally and in its turn
defines the Poisson brackets on the space of smooth functions on
$T^\ast M^n$ or on open domains in $T^\ast M^n$:
\begin{equation}
\{f,g\}=\sum_i
\left(
\frac{\partial f}{\partial x^i}
\frac{\partial g}{\partial p_i} -
\frac{\partial f}{\partial p_i}
\frac{\partial g}{\partial x^i}
\right).
\label{poisson}
\end{equation}

The geodesic flow is a Hamiltonian system on $T^\ast M^n$ with
the Hamiltonian function
$$
H(x,p) = \frac{1}{2} g^{ij}(x)p_i p_j.
$$
This means that the evolution of any function $f$ along trajectories
of the system is given by the Hamiltonian equations
$$
\frac{d f}{dt} = \{f,H\}.
$$
If a function $f$ is preserved by the flow, i.e.
$$
\frac{df}{dt} = \{f,H\} = 0,
$$
it is said that this function is a first integral of the system.

Since the Poisson brackets are skew-symmetric, the function $H$ is
a first integral. This implies that the set of unit momenta vectors
$S M^n$ is invariant under the flow:
$$
S M^n = \{(x,p): |p| = \sqrt{g^{ij}(x)p_i p_j} = 1\} =
\left\{ H=\frac{1}{2} \right\}.
$$
The restrictions of the geodesic flow onto different level sets $H =
{\rm const} \neq 0$ are smoothly trajectory equivalent and this equivalence
is established by constant reparametrization depended only on the values
of $H$. Therefore it is enough to consider the flow only on $S M^n$.

Take the Liouville measure on $S M^n$. This means that the measure of
a set $U \subset S M^n$ is defined as
$$
\mu(U) = \int_x \mu(U \cap S_x) \sqrt{\det g}\,
d x^1 \wedge \dots \wedge d x^n
$$
where $S_x$ is the $(n-1)$-dimensional sphere of unit covectors at the point
$x \in M^n$. In this event the measure on $S_x$
coincides with the measure on the unit sphere in $\R^n$ and this
coincidence is established by an orthogonal map $T^\ast_x M^n \to \R^n$.

B) {\sl Integrability of geodesic flows} \cite{BT,T1}.

The geodesic flow is called {Liouville} integrable if
in addition to $I_n = H$
there are $(n-1)$ first integrals $I_1, \dots, I_{n-1}$ defined on $S M^n$
such that

1) the integrals $I_1,\dots,I_n$ are in involution:
$\{I_j, I_k\} = 0$,

2) the integrals $I_1,\dots,I_{n-1}$ are functionally independent on the
full measure subset $W \subset S M^n$.

To define the Poisson brackets correctly we
extend $I_1, \dots, I_{n-1}$ onto
a neighborhood of $S M^n \subset T^{\ast}M^n$ as follows
$$
I_j (x,p) = I_j (x,p/|p|), \ \ \ j=1,\dots,n-1.
$$
Vanishing of the Poisson brackets of these functions on $S M^n$
does not depend on the choice of $f$.

If the metric and the first integrals $I_1,\dots,I_{n-1}$ are
real-analytic we say that the flow is analytically integrable.

If the geodesic flow is integrable, then a full measure subset $\widetilde{W}$
of $W \subset S M^n$ is foliated by invariant $n$-dimensional tori and
moreover for any such a torus there is its neighborhood
$U \subset W \subset S M^n$ such that

1) there are coordinates
$\varphi_1,\dots,\varphi_n$ defined modulo $\Z$ and $I_1,\dots,I_{n-1}$ in
$U$;

2) every level set $\{I_1 = c_1,\dots, I_{n-1} = c_{n-1}\}$ is an invariant
(Liouville) torus;

3) the flow is linearized in these coordinates as follows:
\begin{equation}
\dot{\varphi}_1 = \omega_1(I_1,\dots,I_{n-1}),\ \dots \ ,
\dot{\varphi}_n = \omega_n(I_1,\dots,I_{n-1}),
\label{linear}
\end{equation}
$$
\ I_1 = {\rm const}, \ \dots\ , I_{n-1} = {\rm const}.
$$
This subset $\widetilde{W}$ is distinguished as the preimage of the set
of regular values of the momentum map $SM^n \to \R^{n-1}$:
$$
x \to (I_1(x),\dots,I_{n-1}(x)).
$$

C) {\sl Entropy.}

Let $X$ be a compact space
and $T: X \to X$ be its homeomorphism.

Take an invariant Borel measure $\mu$ on $X$ such that $\mu(X) < \infty$.
For any a disjoint measurable countable decomposition
$$
X = \sqcup\, U_i
$$
the entropy of the decomposition is defined by the following formula
$$
h(U) = - \sum \mu(U_i) \log \mu(U_i)
$$
assuming that $\mu(U_j) \log \mu(U_j) = 0$ for $\mu(U_j) = 0$.
Let $\{U_i\}$ be such a decomposition.
For any $k \in \N$ define the decomposition $\wedge^k U$ as follows:
$$
X = \sqcup\, U_{i_0 \dots i_{k-1}}
$$
where
$$
x \in U_{i_0 \dots i_{k-1}} \ \ \ \mbox{iff} \ \ \
x \in U_{i_0}, Tx \in U_{i_1},
\dots, T^{k-1}x \in U_{i_{k-1}}.
$$
Now put
$$
h_{\mu}(U,T) = \limsup_{k \to \infty} \frac{h(\wedge^k U)}{k}
$$
and define the measure entropy of $T$ with respect to $\mu$
(the Kolmogorov--Sinai entropy) as
$$
h_{\mu}(T) = \sup_{U \ \mbox{\small with} \
h_{\mu}(U,T) < \infty} h_{\mu}(U,T).
$$

To any open covering
$$
X \subset \cup\, V_j
$$
of $X$
corresponds the series of coverings
$\wedge^k V$ defined as follows:
$$
X \subset \cup\, V_{j_0 \dots j_{k-1}}
$$
where
$$
x \in V_{j_0 \dots j_{k-1}}\ \ \ \mbox{iff} \ \ \
x \in V_{j_0}, Tx \in V_{j_1},
\dots, T^{k-1}x \in U_{j_{k-1}}.
$$
Usually $\wedge^k V$ contains subsets which still form coverings of $X$
and for any $k \in \N$ put $C(k,V,T)$ to be a minimal cardinality of such
a subset. Now put
$$
h(V,T) = \limsup_{k \to \infty}\frac{\log C(k,V,T)}{k}
$$
and define the topological entropy of $T$ as
$$
h_{\rm top}(T) = \sup_{V} h(V,T).
$$

By the Bowen theorem, $h_{\rm top}(T)$ equals the supremum of the
measure entropies with respect to invariant ergodic Borel measures
$\mu$ such that $\mu(X) = 1$.

{\sl Example.} Let $A$ be an automorphism of a torus $T^n = \R^n/\Z^n$
given by a matrix $A \in SL(n,\Z)$. Take a coordinates $x^1,\dots,x^n$
on $T^n$ such that these coordinates are defined modulo $\Z$,
the automorphism $A$ is linear in terms of $x^1,\dots,x^n$ and
$$
\int_{T^n} dx^1 \wedge \dots \wedge dx^n =1.
$$
Then the topological entropy of $A$ and the measure entropy with respect to
$d\mu = dx^1 \wedge \dots \wedge dx^n$ coincide and equal
$$
h_{\rm top}(A) = h_{\mu}(A) = \log \left(
\max_{\lambda \in \Sp A} |\lambda|\right).
$$
Therefore, $h_{\rm top}(A)$ vanishes if and only if all eigenvalues of $A$
lies on the unit circle in $\C$.

D) {\sl The entropies of geodesic flows.}

Let
$$
F_t: S M^n \to S M^n
$$
be a translation along trajectories per the time $t$.

By the definition, the entropy of the geodesic flow is the entropy of the
map
$$
T: S M^n \to S M^n
$$
which is the translation along trajectories per unit time: $T = F_1$.

Recall the definition of Lyapunov exponents. Let $v$ be a tangent vector to
$S M^n$. For any such a vector its norm $|v|$ is defined as follows. Let
$v \in T_q S M^n$ and decompose it into the sum $v = v_M + v_S$,
where $v_M$ is the component tangent to $M^n$ and $v_S$ is the component
tangent to $S_x$ where $q = (x,p) \in S M^n$.
 As in the definition of the Liouville measure, $S_x$
is endowed with a metric by an orthogonal map $T_x M^n \to \R^n$.
Now put
$$
|v|^2 = |v_M|^2 + |v_S|^2
$$
where norms of $v_M$ and $v_S$ are defined by the metrics on $M^n$ and
$S_x$.

On the full measure subset $U$ of $S M^n$
there is a correctly defined map from nonzero tangent
vectors at the points of $U$ to $\R$:
$$
v \longrightarrow \limsup_{t \to \infty} \frac{\log |F_t^\ast (v)|}{t}.
$$
At any point $q \in U \subset S M^n$ such a map takes $2n-1$ values
$$
l_1 \leq l_2 \leq \dots \leq l_k  \leq 0 \leq l_{k+1} \leq \dots \leq
l_{2n-2}
$$
where the zero value is attained on the vector tangent to
the trajectory of the flow. Another values $l_1,\dots,l_{2n-2}$ are
called Lyapunov exponents and some of them may coincide with another.
The number of negative Lyapunov exponents depends on $q$.

The Pesin formula for the measure entropy of the geodesic flow with
respect to any smooth invariant measure $\mu$ on $S M^n$ reads
$$
h_{\mu} = - \int_{S M^n} \sum_{j=1}^{k(q)} l_j(q) d\mu.
$$
It is evident that for the flow
(\ref{linear}) its Lyapunov exponents vanish. Since an
integrable geodesic flow
has such a behavior on a full measure set, the Pesin formula implies
that the entropy of an integrable
flow vanish for any smooth invariant measure on
$S M^n$ and, in particular, for the Liouville measure.

This already follows from the inequality
$$
h_{\mu} \leq - \int_{S M^n} \sum_{j=1}^{k(q)} l_j(q) d\mu.
$$
first established by Margulis in the middle of the 1960s.

\section{The construction of the metric and the lower
estimate for the entropy}

The construction of the metric on $M_A$ is as follows.

Take linear coordinates $x^1,\dots,x^n$ on $T^n$ for which
the map $A$ is linear and take a coordinate $z$ on $\R/\Z$.
These are coordinates on an infinite cylinder ${\cal C} = T^n \times \R$
which descend to coordinates on $M_A$,
the quotient of ${\cal C}$ with respect to the $\Z$-action
generated by
\begin{equation}
(X,z)  \to (AX,z), \ \ \ X = (x^1,\dots,x^n)^\top.
\label{maction}
\end{equation}
The symplectic form takes the form
\begin{equation}
\omega  = \sum_{i=1}^n d x^i \wedge d p_i + dz \wedge dp_z.
\label{variables}
\end{equation}

Define the metric
$$
d s^2 = g_{jk}(z) d x^j d x^k + d z^2
$$
where
\begin{equation}
G(z) = (g_{jk}(z)) =
\gamma(z)^\top \widehat{G} \gamma(z).
\label{metric}
\end{equation}
where $\widehat{G}$ is an arbitrary positive
symmetric $n \times n$-matrix and $\gamma(z)$ is an analytic curve in
$SL(n,\R)$ satisfying the two following properties:
$$
\gamma(z+1)=\gamma(z) A^{-1} \quad \mbox{and} \quad \gamma(0)=E.
$$

It is easily seen that such a curve always exists. Indeed, if all the
eigenvalues of $A$ are positive, then it suffices just to put
$\gamma(z)=e^{-zG_0}$, where $e^{G_0}=A$. 

If the matrix $G_0=\log A \in sl(n,\R)$ does not exist,
then we can use the following simple
construction. Decompose $A$ into product 
of matrices $A_1$ and $A_2$ such that

1) $A=A_1 A_2$;

2) there are $G_i \in sl(n,\R)$ such that $e^{G_i}=A_i$,
$i=1,2$;

3) $A_2$ commute with $e^{zG_1}$ for any $z$ (in particular,
$A_1$ and $A_2$ commute).

To prove that such a decomposition exists take a Jordan form of $A$,
which is a block matrix. Take now a diagonal matrix $A_2$, whose
entires equal $\pm 1$ and such that all eigenvalues of $A A_2 = A A_2^{-1}$
has positive eigenvalues. 
Since $\det A=1$, the matrix $A_2$ has an even number of diagonal elements,
which equal $-1$, and therefore there is a matrix 
$G_2 \in so(n) \subset sl(n,\R)$
such that $A_2 = e^{G_2}$.
Now it remains to put $A_1 = A A_2^{-1}$.

Given $A_1$ and $A_2$, put $\gamma(z)=e^{-zG_2}e^{-zG_1}$.

It is clear that (\ref{metric}) defines a metric on an infinite cylinder
${\cal C}$ and the metric is invariant with respect
to the action (\ref{maction}).
Therefore, this metric descends to a metric on the quotient space
$M_A = {\cal C}/{\Z}$.

\begin{lemma}l
The geodesic flow of the metric
(\ref{metric}) on the cylinder ${\cal C}$ is integrable, i.e., it admits
$n+1$ first integrals
$$
I_1 = p_1, \ \dots, \ I_n = p_n, \
I_{n+1} = H =
\frac{1}{2} \left(g^{ij}(z)p_i p_j + p_z^2 \right)
$$
which are in involution and for any open subset $U \subset T^\ast {\cal C}$
these integrals are functionally independent on a full measure subset of $U$
with respect to the Liouville measure.
\end{lemma}

{\sl Proof.}
It is clear that these integrals are functionally independent at least
on the set where $p_z \neq 0$. By (\ref{poisson}) and (\ref{variables}),
the momenta variables
are in involution:
$$
\{ p_i, p_j \} = 0, \ \ i,j=1,\dots,n,
$$
and, moreover, since $H$ does not depend on $x^1,\dots,x^n$,
we have
$$
\{p_i, H\} = 0, \ \ i=1,\dots,n.
$$
This proves the lemma.

Now take a torus $T^n \subset S M^n$ formed by the points with $z=0$ and
$p_1 = \dots = p_n = 0$. Since $p_1,\dots,p_n$ are preserved by the flow,
the translation $T =F_1$ along trajectories of the geodesic flow
per unit time maps its torus into itself:
$$
(X,0) \to (X,1) \sim (AX,0)
$$
and we see that the dynamical system $T: S M_A \to S M_A$ contains a
subsystem isomorphic to the torus automorphism $A: T^n \to T^n$.
It is known that the topological entropy of a system is not less than the
topological entropy of any of its subsystems. Therefore we conclude that
$$
h_{\rm top}(T) \geq h_{\rm top}(A) = \log \left(
\max_{\lambda \in \Sp A} |\lambda| \right).
$$

For proving integrability of the flow
we are left to descend the first integrals $p_1,\dots,p_n$
to $S M_A$.
We can not do that straightforwardly but may substitute them
by some functions of $p_1,\dots,p_n$  which are invariant under
the action of $A$ and functionally
independent almost everywhere.

\section{Proof of Theorem 1 for $A$ with real eigenvalues}

The action of $A$ on $M_A$ generates the natural action on tangent vectors,
the differential. We expand the action of $A$ onto $T^\ast M_A$ by
assuming that $A$ preserves the form $\omega$. This action is also linear in
terms of $p_1,\dots,p_n$. Denote this action by $\widetilde{A}$.
It is uniquely define by the equation
$$
\left(
\begin{array}{cc}
A^\top & 0 \\
0 & \widetilde{A}^\top
\end{array}
\right)
\left(
\begin{array}{cc}
0 & 1 \\
-1 & 0
\end{array}
\right)
\left(
\begin{array}{cc}
A & 0 \\
0 & \widetilde{A}
\end{array}
\right) =
\left(
\begin{array}{cc}
0 & 1 \\
-1 & 0
\end{array}
\right)
$$
which means that $\omega$ is preserved and reads
$$
A^\top \widetilde{A} = 1.
$$

Let all eigenvalues of $A$ be real.
Then all eigenvalues of $\widetilde{A}$ are real.
Take linear coordinates
$p_1,\dots,p_n$ such  that $\widetilde{A}$ attains its Jordan form:
$$
\widetilde{A} =
\left(
\begin{array}{cccc}
B_0 & 0 & \dots & 0 \\
0 & B_1 & \dots & 0 \\
0 & \dots & \dots & 0 \\
0 & \dots & 0 & B_k
\end{array}
\right)
$$
where
$B_0$ is a diagonal matrix
$$
B_0 = {\rm diag}(\mu_1,\dots,\mu_l)
$$
and for $j \geq 1$
each matrix $B_j$ is an $n_j \times n_j$-matrix of the form
$$
B_j = \left(
\begin{array}{cccccc}
\lambda_j & 1 & 0 & \dots & 0 & 0 \\
0 & \lambda_j & 1 & \dots & 0 & 0 \\
& & & \dots & & \\
0 & 0 & 0 & \dots & \lambda_j & 1 \\
0 & 0 & 0 & \dots & 0 & \lambda_j
\end{array}
\right)
$$
where $n_1+\dots + n_k + l =n$.
Hence redenote the variables as follows
$$
p_1,\dots,p_n
\longrightarrow
q_1,\dots,q_l,p_{11},\dots,p_{1n_1},\dots,p_{k1},
\dots,p_{kn_k}.
$$

Introduce the following polynomial
$$
Q = q_1 \dots q_l p_{11}^{n_1} \dots p_{k1}^{n_k}.
$$
Since $A \in SL(n,\Z)$, we have $A^\top \in SL(n,\Z)$ and, therefore,
$\widetilde{A} = (A^\top)^{-1} \in SL(n,\Z)$. This implies
$$
\det \widetilde{A} = \mu_1 \dots \mu_l \lambda_1^{n_1} \dots \lambda_k^{n_k}
=1.
$$
Since
$$
Q \to (\mu_1 q_1) \dots (\mu_l q_l) (\lambda_1 p_{11})^{n_1} \dots
(\lambda_k k_{k1})^{n_k}
= (\mu_1 \dots \mu_l \lambda_1^{n_1} \dots \lambda_k^{n_k}) Q
$$
this results in the following lemma.

\begin{lemma}
The polynomial $Q$ is an invariant of the action $\widetilde{A}$.
\end{lemma}

Before constructing the full family of first integrals let us
prove the technical lemma which we shall need.

\begin{lemma}
\label{lemmanil}
Let $L$ be an operator acting on the ring
$\R[p_1,\dots,p_n]$
of polynomials in $p_1,\dots,p_n$ as follows:
\begin{equation}
L \cdot f(p_1,\dots,p_n) = f(L \cdot p_1,\dots,
L \cdot p_n), \ \ \ \ f \in \R[p_1,\dots,p_n],
\label{action}
\end{equation}
where
\begin{equation}
L \cdot p_1 = \lambda p_1, \ \ L \cdot p_k = \lambda p_k + p_{k-1}
\ \ \mbox{for $k=2,\dots,n$}
\label{Aaction}
\end{equation}
and $\lambda$ is a constant.

Then for any $k=1,\dots,n$ there is a polynomial $G_k \in \R[p_1,\dots,p_n]$
of degree $k$ such that

1) $G_k$ depends only on $p_1,\dots,p_{k+1}$ and has the form
$$
p_{k+1} H_{k1}(p_1,\dots,p_k) + H_{k2}(p_1,\dots,p_k)
$$
where $H_{k1}, H_{k2} \in \R[p_1,\dots,p_k]$;

2) the operator $L$ acts on $G_k$ as follows
$$
L \cdot G_k = \lambda^k G_k + p_1^k.
$$
\end{lemma}

{\sl Proof.}
Let $V^l_k$ be the space of homogeneous polynomials in $p_1,\dots,p_l$ of
degree $k$. It is clear from (\ref{action}) and (\ref{Aaction})
that $L(V^l_k) \subset V^l_k$.

Notice that the linear operator
\begin{equation}
(L - \lambda^k): V^l_k \to V^l_k
\label{action-k}
\end{equation}
is nilpotent.
Indeed, let us introduce the following order on monomials from $V^l_k$:
$$
p_1^{\alpha_1} \dots p_l^{\alpha_l} \prec p_1^{\beta_1} \dots p_l^{\beta_l}
\ \ \
\mbox{if}
\ \ \
\alpha_r = \beta_r \ \ \mbox{for $l > m$ and} \ \ \alpha_m < \beta_m.
$$
Then $L$ acts on any monomial $F = p_1^{\alpha_1} \dots p_k^{\alpha_k}$
as follows
$$
L \cdot F = \lambda^k F + \sum_j D_j
$$
where $D_j$ are monomials such that $D_j \prec F$.

It is also clear that the kernel of the action (\ref{action-k})
is generated by $p_1^k$.

This implies that in some basis $e_1,\dots,e_N$
for $V^l_k$ $L$ takes the form
\begin{equation}
\left(
\begin{array}{cccccc}
\lambda^k & 1 & 0 & \dots & 0 & 0 \\
0 & \lambda^k & 1 & \dots & 0 & 0 \\
& &  & \dots & & \\
0 & 0 & 0 & \dots & \lambda^k & 1 \\
0 & 0 & 0 & \dots & 0 &\lambda^k
\end{array}
\right),
\label{form}
\end{equation}
where $e_1 = p_1^k$.

Put $F_k = p_{k+1} p_1^{k-1}$.
Then we have
$$
L \cdot F_k = \lambda^k F_k + \lambda^{k-1}p_k p_1^{k-1}.
$$
Look for solutions $H_k$ and $c_k$ to the equation
\begin{equation}
(L - \lambda^k) \cdot H_k = c_k p_1^k - \lambda^{k-1} p_k p_1^{k-1}
\label{G}
\end{equation}
where $H \in V^k_k$ and $c_k \in \R$.
In some basis $e_1,\dots,e_N$ for $V_k^k$ the operator $L$ has the form
(\ref{form}) and,
since the monomial $ p_k p_1^{k-1}$ is not maximal in $V_k^k$,
$$
p_k p_1^{k-1} = \sum_{j \leq (N-1)} a_j e_j.
$$
The vectors $e_2,\dots,e_{N-1}$ lies in the image of $(L -\lambda^k)$
and therefore the equation (\ref{G}) is solvable in $H_k$ for
$c_k = a_1 \lambda^{k-1}$. Take a solution $H_k$ to it.
We see that $F_k+H_k = p_{k+1} p_1^{k-1}$ satisfy the equation
$$
L \cdot (F_k + H_k) = \lambda^k (F_k + H_k) + c_k p_1^k.
$$
If $C_k = 0$ then $(F_k+H_k)$ lies in the kernel of $(L - \lambda^k)$ but
$(F_k+H_k)$ is not proportional to $p_1^k$. Hence $c_k \neq 0$ and
we are left to put
$$
G_k = \frac{1}{c_k}\left( F_k + H_k \right).
$$
This proves the lemma.

These are some simplest examples of the polynomials $G_k$:
$$
G_1 = p_2, \ \ \ G_2 = p_2^2 - 2p_1 p_3, \ \ \
G_3 = p_2^3 + 3p_1^2 p_4 - 3 p_1 p_2 p_3,
$$
$$
G_4 = p_2^4 - 4 p_1^3 p_5 - 4 p_1 p_2^2 p_3 + 2 p_1^2 p_3^2 + 4 p_1^2 p_2 p_4.
$$

\begin{corollary}
Given an  action $L$ on $\R[p_1,\dots,p_n]$  such that it
satisfies (\ref{action}), this action admits $(n-1)$ almost invariant
rational functions functionally independent outside an algebraic subvariety
of positive codimension. These are
$$
J_k = \frac{G_k}{p_1^k}, \ \ \ k=1,\dots,n-1
$$
which are transformed by $A$ as follows
$$
J_k \stackrel{L}{\longrightarrow} J_k + \frac{1}{\lambda^k}.
$$
\end{corollary}

Notice that the functional independence statement
follows from the fact that
each polynomial $G_k$ depends only on $p_1,\dots,p_{k+1}$ and is
linear in $p_{k+1}$.

Now we are ready to finish the proof of Theorem 1 for $A$ with real
eigenvalues.

Put
$$
I_1 = \exp\left(-\frac{1}{Q^2}\right)
\sin \left(2\pi \frac{\log q_1}{\log \mu_1}\right),\ \dots,
I_l = \exp\left(-\frac{1}{Q^2}\right)
\sin \left(2\pi \frac{\log q_l}{\log \mu_l}\right).
$$

To each series of variables $p_{j1},\dots,p_{jn_j}$
we apply Lemma \ref{lemmanil} and construct the polynomials
$G_1,\dots,G_{n_j-1}$.
Now put
$$
I_{j1} =\exp\left(-\frac{1}{Q^2}\right)
\sin \left(2\pi \frac{\log p_{j1}}{\log \lambda_j}\right), \ \
I_{j2} = \exp\left(-\frac{1}{Q^2}\right)
\sin \left(2\pi \lambda_j \frac{G_1}{p_{j1}} \right),
$$
$$
\dots, \ \ I_{jm} = \exp\left(-\frac{1}{Q^2}\right)
\sin \left(2\pi \lambda_j^{m-1} \frac{G_{m-1}}{p_{j1}^{m-1}} \right),
\ \ \dots,
$$
$$
\dots, \ \
I_{jn_j} = \exp\left(-\frac{1}{Q^2}\right)
\sin \left(2\pi \lambda_j^{n_j-1} \frac{G_{n_j-1}}{p_{j1}^{n_j-1}} \right).
$$
These functions are smooth, invariant under the action of $\widetilde{A}$
and functionally independent at any fiber $S_x {\cal C}$ outside an
algebraic subset of positive codimension. In fact, outside this singular set
where they functionally dependent these functions substitute $p_{j1},\dots,
p_{jn_j}$.

The functions $I_1,\dots,I_l,I_{11},\dots,I_{kn_k}$
are  functionally independent at any fiber $S_x {\cal C}$,  invariant
under $\widetilde{A}$ and,  therefore, descend to functions on
$S M_A$. Since these functions
depend only on the momenta variables, they
are in involution and are first integrals of the geodesic flow on $M_A$.

We conclude that this family gives us a complete family of first integrals
and therefore the geodesic flow on $M_A$ is integrable.

The case of Theorem 1 concerning automorphisms
$A$ with real eigenvalues is established.

\section{Proof of Theorem 1 for $n=2$}

The case when all eigenvalues are real is already considered.
In fact, the case when $A$ is not diagonalized and therefore
in a convenient coordinates equals
$$
\left(
\begin{array}{cc}
1 & 1 \\
0 & 1
\end{array}
\right)
$$
was the initial one discovered by Butler \cite{Butler}
and the case when $A$ is diagonalized with real eigenvalues was
considered by us in \cite{BT}.

Hence we assume that $\lambda$ and $\bar{\lambda}$ are eigenvalues of $A$
and, since $A \in SL(2,\Z)$, we have
$$
\lambda + \bar{\lambda} \in \Z, \ \ |\lambda| =1.
$$
This means that $\lambda = \cos \varphi + i \sin \varphi$ and
$2 \cos \varphi \in \Z$. The latter inclusion implies $\cos \varphi \in
\{ \pm 1, \pm 1/2, 0\}$. If $\cos\varphi = \pm 1$ then $\lambda = \pm 1$
and hence $\lambda$ is real. Therefore we are left with the
following cases: in the momenta coordinates $p_1,p_2$ the action
$A$ is a rotation by
$$
\varphi =
\pm \frac{\pi}{2}, \ \pm \frac{2\pi}{3}, \ \pm \frac{\pi}{3}.
$$
It is clear that this action preserves
$$
I_1(p_1,p_2) = p_1^2 + p_2^2.
$$
Put
$$
\psi = \arcsin \frac{p_2}{\sqrt{p_1^2 + p_2^2}}
$$
and notice that $A$ acts as
$$
\psi \to \psi + \varphi.
$$
Now we put
$$
I_2(p_1,p_2) = {\rm Re}\, (p_1 + ip_2)^k,
$$
where $\varphi = \pm 2\pi/k$.

It is easy to notice that the functions $I_1$ and $I_2$ are functionally
independent almost everywhere.

This proves Theorem 1 for $n=2$.

\section{Proof of Theorem 2}

Take linear coordinates $u$ and $v$ on $T^2$
such that $A$ of the form (\ref{matrix}) acts as
\begin{equation}
u \to \lambda^{-1} u, \ \ v \to \lambda v
\label{action2}
\end{equation}
with
$$
\lambda = \frac{3+\sqrt{5}}{2},
$$
and also take a linear coordinate $z$ on $S^1$ which is lifted to
a coordinate on $M_A$ defined modulo $\Z$. These coordinates are
completed by $(p_u,p_v,p_z)$ to coordinates on $T^{\ast}M_A$ such that
the symplectic form on the cotangent bundle is
$$
\omega = du \wedge dp_u + dv \wedge dp_v + d z \wedge dp_z
$$
and $A$ acts on the momenta as
\begin{equation}
p_u \to \lambda p_u, \ \ p_v \to \lambda^{-1}p_v, \ \ p_z \to p_z.
\label{momenta}
\end{equation}
Now the metric on $M_A$ is
$$
ds^2 = dz^2 + e^{2z \log \lambda} du^2 + e^{-2z\log \lambda} dv^2
$$
and the Hamiltonian function on $T^{\ast}M_A$ is
$$
H = \frac{1}{2} \left( p_z^2 + e^{-2z \log \lambda} p_u^2 +
e^{2z\log \lambda} p_v^2 \right).
$$
There are three functionally independent almost everywhere
first integrals of the geodesic flow on the universal covering:
$$
I_1 = p_u, \ \ I_2 = p_v, \ \ I_3 = H.
$$

Consider the restriction of the flow on the compact level surface
$S M_A$ defined as
$$
S M_A = \left\{H =\frac{1}{2}\right\}.
$$

1) If $p_u p_v \neq 0$, then the lift of a trajectory on the
universal covering is trapped in the layer
$$
c_1 e^{-2z \log \lambda} + c_2 e^{2z\log \lambda} \leq 1
$$
with the constants $c_1 = p_u^2$ and $c_2 = p_v^2$.
This layer is invariant under $\Z^2$ actions under translations by
vectors of the lattice $\Lambda$. Here $T^2 = \R/\Lambda$ and
in the coordinates $u$ and $v$ the vectors from $\Lambda$ have
irrational coefficients. There are two different kinds of such
trajectories:

1a) A trajectory,
for which
$$
p_z \neq 0 \ \ \ \mbox{or} \ \ \ p_u^2 \neq p_v^2,
$$
lies on an invariant torus in
$S M_A$ and its Lyapunov exponent vanish.
These inequalities describe the set on which the first
integrals $I_1, I_2$, and $I_s$ are functionally independent;

1b) Trajectories with
$$
p_z = p_u^2 - p_v^2 =0
$$
form a submanifold which is evidently diffeomorphic to two copies of
$M_A$ corresponding to two possibilities: $p_u = \pm p_v$.
Each of this copies is fibered over $S^1$ and this fibration is induced
by (\ref{fibration}).
Since $\dot{z} = p_z = 0$ on such a trajectory, it lies on the level
$z = {\rm const}$ which is a torus with linear coordinates $u$ and
$v$. The flow is linear in these coordinates and has constant velocities.
Therefore the Lyapunov exponents for such a trajectory are zero.

2) Trajectories with $p_v =0$ form a submanifold $N^u$. Since $M_A$ is
parallelizable, we see that $N^u$ is diffeomorphic to $M_A \times S^1$
and the flow on it is described by the equations
\begin{equation}
\dot{p}_u = 0, \ \
\dot{p}_z = \log{\lambda}e^{-2z\log \lambda} p_u^2, \ \
\dot{u} = e^{2z\log \lambda}p_u, \ \
\dot{z} = p_z.
\label{flow}
\end{equation}
There are two invariant submanifolds of $N^u$,
which are
$$
V^+ = \{ p_u = p_v = 0, p_z =1\}, \ \ \
V^- = \{ p_u = p_v = 0, p_z =-1\}.
$$
Any trajectory with $p_v = 0$ satisfies the inequality
$$
e^{-2z \log \lambda} \leq \frac{2}{p_u^2}
$$
and we see that the lift of such a trajectory onto the universal covering
is not trapped into any layer but just bounded in $z$ from below.
Hence

{\sl any trajectory on $S M_A$ with $p_u \neq 0$ and $p_v = 0$ is
asymptotic to a trajectory from $V^+$ as
$t \to \infty$ and asymptotic to a trajectory from $V^-$
as $t \to -\infty$.}

Since the metric is invariant with respect to the
$A$-action and the action of (\ref{flow}) on the tangent vector field
$$
\xi = \frac{\partial}{\partial p_u}
$$
is trivial: $F_t^{\ast}(\xi) = \xi$, we derive from (\ref{action2})
that the Lyapunov exponent corresponding to this vector is positive:
$$
\limsup_{t \to \infty}\frac{\log |F_t^{\ast}(\xi)|}{|\xi|} > 0.
$$

3) The submanifold $N^v$ of $S M_A$
is defined by the equation $p_u = 0$. It is analyzed in completely the
similar manner as $N^u$ and we derive that

{\sl any trajectory from $N^v$
with $p_v \neq 0$ is
asymptotic to a trajectory from $V^-$ as
$t \to \infty$ and asymptotic to a trajectory from $V^+$
as $t \to -\infty$.}

We see that all trajectories in $\{N^u \cup N^v\} \setminus \{V^+ \cup V^-\}$
are not closed which implies that all invariant Borel measures
on $N^u$ and $N^v$ are supported by $V^+ \cup V^-$.  Otherwise it would
contradict to the Katok theorem \cite{Katok}, which reads that given
a compact manifold with an invariant Borel
measure with nonzero Lyapunov exponents
the support of the measure lies in the closure of periodic trajectories.

This finishes the proof of Theorem 2.

There is a natural invariant measure on $V^+$, which is
\begin{equation}
d\mu = du \wedge dv \wedge dz,
\label{measure}
\end{equation}
and the measure entropy with respect to $d\mu$ equals
the topological entropy of the automorphism $A$ of the torus,
which is $\log \lambda$.

By the Bowen theorem, the topological entropy of a flow
equals the supremum of the measure entropies of the flow taken over
all invariant ergodic Borel measures. For an integrable flow with first
integrals $I_1,\dots, I_n$ it is easy to derive from this ergodicity
restriction for measures that there are constants
$C_1,\dots,C_n$ such that this supremum may be taken over all
measures supported on the level $\{I_1 = C_1, \dots, I_n = C_n\}$
(see, for instance, \cite{T3}). Knowing the first integrals of the geodesic
flow on $M_A$ and the behavior of it trajectories, we see that the topological
entropy of this flow is the supremum of the measure entropies supported by
$V^+$ or $V^-$. But the restrictions of the flow onto these sets the
topological entropy equals $\log\lambda$ and this establishes Corollary 1.

In fact, Theorem 2 describes the geodesic flow on the universal
covering of $M_A$, which is the solvable Lie group ${\rm SOL}$. This
manifold  is a model for one of Thurston's canonical three-geometries.
Asymptotic properties of its geodesic flow
were studied in \cite{Leeb} where some general results on solvable
groups were proved, which imply that the Martin boundary
of ${\rm SOL}$ consists in a single point, and in \cite{Troyanov} where a
rather complex ``horison'' of the group ${\rm SOL}$
defined via the asymptotics of geodesics was described.

Speaking about the geodesic flow on $M_A$ we would like to remind the first
integrals of it, which were found in \cite{BT}:
$$
I_1 = p_u p_v, \ \ \
I_2 = \exp \left(-\frac{1}{p_u^2 p_v^2}\right) \sin\left(2\pi
\frac{\log p_u}{\log \lambda}\right), \ \ \
I_3 = H.
$$
It is easy to check from (\ref{momenta}) that these functions are invariants
of $A$ and therefore descend to $S M_A$. They are the first integrals of the
geodesic flow on $S M_A$ which are functionally independent on a full
measure subset of $S M_A$.

\section{Some remarks and open problems}

The problem of topological obstructions to integrability
was posed by Kozlov who also found the first known obstruction:
he proves that if there is an
analytically integrable geodesic flow on an oriented closed
two-di\-men\-sio\-nal
manifold then this manifold is homeomorphic to the two-sphere $S^2$ or
the two-torus $T^2$ \cite{K1,K2}.
As shown by Kolokol'tsov \cite{Kol} this
also true for geodesic flows on two-ma\-ni\-folds, which are integrable
in terms of smooth first integrals, which are real-analytic functions of
the momenta. But the following problem remains unsolved

\begin{problem}
Can the Kozlov theorem be generalized for $C^\infty$ metrics on
two-ma\-ni\-folds
with geodesic flows integrable in terms of $C^{\infty}$ first integrals ?
\end{problem}

Speaking not about integrability but on existence of metrics
whose geo\-de\-sic flows have zero Liouville entropy we would like to remind
the problem posed by Katok:

\begin{problem}
Does there exist a smooth (at least $C^2$) geodesic flow with zero
Liouville entropy on a two-sphere
with $g \geq 2$ handles ?
Or more general, do there exists such a flow on a closed manifold
admitting negatively curved metric ?
\end{problem}

There is a similar question for mappings which also belongs to Katok.

\begin{problem}
Does there exist a smooth (at least $C^{1+\alpha}$) diffeomorphism $f$ of an
$n$-dimensional torus $T^n$ with $n \geq 3$ such that it induces
an Anosov automorphism $f_{\ast}: \Z^n \to \Z^n$
in homologies (and therefore, its topological entropy is positive)
and its measure entropy with respect to some invariant smooth measure on
$T^n$ vanishes ?
\end{problem}

A generalization of the Kozlov theorem for higher-dimensional manifolds
was found in \cite{T1,T2} where it was shown that if the geodesic flow on
a closed manifold $M^n$ is analytically integrable then the unit cotangent
bundle $S M^n$ contains an invariant torus $T^n$ such that its projection
onto the base
$$
\pi: T^n \subset S M^n \to M^n
$$
induces a homomorphism of the fundamental groups
$\pi_{\ast}: \pi_1(T^n) \to \pi_1(M^n)$
whose image $\pi_{\ast}(\pi_1(T^n))$ has a finite index in $\pi_1(M^n)$:
$$
[\pi_1(M^n) : \pi_{\ast}(\pi_1(T^n))] < \infty.
$$
This implies that

1) the fundamental group of $M^n$ is almost commutative;

2) if the first Betti number $b_1(M^n)$ of $M^n$ equals $k$: $b_1(M^n) = k$,
then the real cohomology ring $H^{\ast}(M^n;\R)$ of $M^n$ contains a subring
isomorphic to the real cohomology ring of the $k$-dimensional torus:
$$
H^{\ast}(T^k;\R) \subset H^{\ast}(M^n;\R).
$$
In particular, this implies that
\begin{equation}
b_1(M^n) \leq n = \dim M^n;
\label{inequal}
\end{equation}

3) if $b_1(M^n) = \dim M^n$, then $H^{\ast}(T^n;\R) = H^{\ast}(M^n;\R)$.

This result is valid for more general case when the flow is not
analytically integrable but so-called geometrically simple and
also is immediately generalized for superintegrable cases when there are more
than $n$ functionally independent real analytic
first integrals and generic tori are
$l$-dimensional with $l <n$ (in this case the ``maximal'' torus whose
fundamental group projects into a group with finite index is $l$-dimensional).

As shown by Butler \cite{Butler} some of
these topological properties do not obstruct $C^{\infty}$ integrability:
for Butler's manifold we have $b_1=2$ and the fundamental group is not almost
commutative and $H^{\ast}$ contains no subring isomorphic to $H^{\ast}(T^2;\R)$
but the inequality (\ref{inequal}) is valid. In fact this is true also for
the geodesic flows on $M_A$ where $A$ is not of finite order.

We would like to introduce the following

\begin{conjecture}
Let the geodesic flow on a Riemannian manifold $M^n$ is integrable in terms of
$C^{\infty}$ first integrals. Then the inequalities
\begin{equation}
b_k (M^n) \leq b_k (T^n) = \frac{n!}{k!(n-k)!},
\label{inequal2}
\end{equation}
hold.
\end{conjecture}

These inequalities mean that
homologically $M^n$ is dominated by the $n$-dimensional
torus. They were already mentioned in talks of the second author (I.A.T.)
in the early 90s. It was derived by Paternain from results of Gromov and
Yomdin that if the topological entropy of
the geodesic flow of a $C^{\infty}$ metric on a simply connected manifold
vanishes, then this manifold is rationally elliptic (in the sense of
Sullivan) \cite{P1} and he also  mentioned that, by results
Friedlander and Halperin, rational ellipiticity implies the inequalities
(\ref{inequal2}).

Actually, it was Paternain who proposed the entropy approach to
finding topological obstructions to integrability. He proposed to
split this problem into two ones: proving the vanishing of the
topological entropy of an integrable geodesic flow and finding
topological obstructions to vanishing of the topological entropy of a flow.
The second problem was already studied and in addition to the results of
Gromov and Yomdin, which we already mentioned above, we would like
to remind the theorem of Dinaburg who proved that if the fundamental group
of the manifold has an exponential growth, then the topological entropy of
the geodesic flow of any smooth metric on the manifold is positive
\cite{D}.

Paternain found some conditions mainly concerning existence of
rather good action-angle variables on the set, where the first integrals
are functionally dependent, which in addition to integrability
imply the vanishing of the topological entropy \cite{P1,P2}
(after that some other similar conditions were exposed in \cite{T3}).

He also conjectured that the topological entropy of an integrable
geodesic flow vanishes and that the fundamental group of a manifold with
an integrable geodesic flow has a subexponential growth.

In \cite{BT} we disproved both these conjectures in the $C^{\infty}$ case.
Since it is proved in \cite{T1}, that if the geodesic flow is
analytically integrable, then the fundamental group of the manifold
has a polynomial growth, we are left with the following real-analytic
version of Paternain's conjecture:

\begin{conjecture}
If the geodesic flow on a closed manifold is analytically integrable, then
the topological entropy of the flow vanishes.
\end{conjecture}

We already mentioned about eight Thurston's
canonical three-geo\-met\-ri\-es, which are the homogeneous geometries of
$S^3, \R^3, H^3, S^2 \times \R, H^2 \times \R, {\rm NIL}$, ${\rm SOL}$, and
$SL(2,\R)$. Here we denote by $H^n$ the $n$-th dimensional Lobachevsky
space. Since the Lyapunov exponents does not vanish at any point, there
are no compact quotients of $H^3$ and $H^2 \times \R$ with integrable
geodesic flows. There are well-known examples of compact quotients of
$\R^3$ and $S^2 \times \R$ with integrable geodesic flows, which are, for
instance, flat tori $T^3$ and $S^2 \times S^1$. The geodesic flow of the
Killing metric on $SU(2) = S^3$ is also integrable. As shown in \cite{Butler}
and \cite{BT} there are compact quotients of ${\rm NIL}$ and ${\rm SOL}$
with integrable geodesic flows. Hence it remains to answer the following
question:

\begin{problem}
Do there exist compact quotients of $SL(2,\R)$ with integrable geodesic
flows ?
\end{problem}

\vskip1cm

{\sl Acknowledgement.}
The authors were supported by the Russian Foundation of Basic Researches
(grants 96-15-96868 and 98-01-00240 (A. V. B.), and 96-15-96877 and
98-01-00749 (I.A.T.)).

\end{document}